\documentclass[12pt]{article} 
\usepackage[cp1251]{inputenc}
\usepackage[ukrainian]{babel}
\usepackage{latexsym,amsfonts,amssymb,amsmath,epsfig}
\usepackage{graphicx,graphics,hhline,cite}
\usepackage{euscript}

\newtheorem{theorem}{Теорема}

\newtheorem{corollary}{Наслідок}

\newtheorem{lemma}{Лема}
\newtheorem{remark}{Зауваження}

\begin{document}

\noindent \textbf{A.\,S. Serdyuk, I.\,V. Sokolenko}\\ {\small (Institute of Mathematics, National Academy of Sciences of Ukraine, Kyiv)}\\
\noindent{serdyuk@imath.kiev.ua, sokol@imath.kiev.ua}
\vskip 2mm
\begin{center}
	\textbf{Approximation by Fourier sums in classes of Weyl--Nagy differentiable  functions  with high exponent of smoothness
		\footnote{This work was partially supported by the Project ''Innovative methods in the theory of differential equations, computational mathematics and mathematical modeling''. 
		}
	}
\end{center}
\vskip 2mm
{\small We establish asymptotic estimates for the least upper bounds
	of approximations in the uniform metric by Fourier sums of order $n-1$ of classes of $2\pi$-periodic Weyl--Nagy differentiable functions, $W^r_{\beta,p}, 1\le p\le \infty, \beta\in\mathbb{R},$ for high exponents of smoothness $r\ (r-1\ge \sqrt{n})$.  We obtain similar estimates in metrics of
	the spaces $L_p, 1\le p\le\infty,$ for functional classes  $W^r_{\beta,1}$.   
}
\vskip 2mm

\noindent \textbf{А.\,С. Сердюк, І.\,В. Соколенко} {\small (Iнститут математики НАН України, Київ)} 
\begin{center}
\textbf{Наближення сумами Фур'є на класах диференційовних в сенсі Вейля--Надя функцій з високим показником гладкості 
}\end{center}

\vskip 2mm
{\small Встановлено асимптотичні оцінки точних верхніх меж відхилень в рівномірній метриці частинних сум Фур'є порядку $n-1$ на класах $2\pi$-періодичних функцій, диференційовних в сенсі Вейля--Надя, $W^r_{\beta,p}, 1\le p\le \infty, \beta\in\mathbb{R},$ при високих показниках гладкості $r\ (r-1\ge \sqrt{n})$. Аналогічні оцінки встановлено і в метриках просторів $L_p, 1\le p\le\infty,$ для функціональних класів $W^r_{\beta,1}$.} 
\vskip 5mm



\vskip 5mm

\textbf{1. Вступ.} Нехай $L_p, \ 1\le p<\infty,$ --- простір
$2\pi$-періодичних сумовних в $p$-му степені на $[-\pi, \pi)$ функцій $\varphi$ зі стандартною нормою
$$
\|\varphi\|_{p}=\bigg(\int\limits_{-\pi}^{\pi}|\varphi(t)|^p dt \bigg)^{1/p};
$$
$L_{\infty}$ ---
простір $2\pi$-періодичних вимірних та істотно обмежених функцій
$\varphi$ в якому  норма задана рівністю
$$ \|\varphi\|_{\infty}=\mathop {\rm ess
	\sup}\limits_{t} |\varphi(t)|;
$$
$C$ --- простір $2\pi$-періодичних неперервних функцій $\varphi$ з нормою
$$
\|\varphi\|_{C}=\max\limits_{t}|\varphi(t)|. $$

Нехай, далі, $W^r_{\beta,p}, r>0, \beta\in\mathbb{R}, 1\le p\le \infty,$
--- класи $2\pi$-періодичних функцій $f$, що зображуються у вигляді згортки
\begin{equation}\label{1}
f(x)=\frac{a_0}{2}+\left(\varphi\ast B_{r,\beta}\right)(x)=\frac{a_0}{2}+\frac{1}{\pi}\int\limits_{-\pi}^{\pi}
\varphi(x-t) B_{r,\beta}(t)dt, \ \ \ a_0\in\mathbb R, 
\end{equation}
з ядрами Вейля-Надя $B_{r,\beta}(\cdot)$ вигляду
\begin{equation}\label{2}
B_{r,\beta}(t)=\sum\limits_{k=1}^\infty k^{-r}\cos\left(kt-\frac{\beta\pi}2\right),\quad r>0,\quad \beta\in\mathbb{R},
\end{equation}
функцій $\varphi$, що задовольняють умову $\varphi\in B_p^0$, де
\begin{equation}\label{3}
 B_p^0:=\left\{\varphi\in L_p: \|\varphi\|_p\le1,\   \int\limits_{-\pi}^{\pi}\varphi(t) dt=0\right\}.
\end{equation}
Класи $W^r_{\beta,p}$  називають класами Вейля--Надя (див., наприклад, \cite{Stepanets1995, Stepanets2005, Sz.-Nagy1938, Stechkin1956}), а функцію $\varphi$ в зображенні \eqref{1} називають $(r,\beta)$-похідною в сенсі Вейля--Надя функції $f$ і позначають через $f^r_\beta$.

При всіх $1<p,s<\infty$, $r>(\frac1p-\frac1s)_+=\left\{\begin{array}{ll}
	\frac1p-\frac1s,& p<s,\\
	0,& p\ge s,
\end{array}\right.$ i $\beta\in\mathbb{R}$ має місце вкладення $W^r_{\beta,p}\subset L_s$ (див., наприклад, \cite[Ch.~V.4]{Stepanets1995}, \cite[Ch.~VI.6]{Stepanets2005}). Окрім того, при довільних $1\le p\le\infty, r>\frac1p, \beta\in\mathbb{R}$ виконуються вкладення $W^r_{\beta,p}\subset C$ (див., наприклад, \cite{Temlyakov1993}).

Якщо  $r\in\mathbb N$ i $ \beta=r,\ $ то   функції вигляду (\ref{2})  є відомими ядрами Бернуллі, а відповідні класи $W^r_{\beta,p}$ збігаються з відомими класами $W^r_{p}$  $2\pi$-періодичних функцій $f$, які мають абсолютно неперервні похідні до $(r-1)$-го порядку включно і такі, що $\|f^{(r)}\|_p\le1$. При цьому 
майже скрізь виконується рівність $f^{(r)}(\cdot)=f^r_\beta(\cdot).$

Для довільної множини $\mathfrak N\subset X$, де $X=C$ або $L_s, \ 1\le s\le\infty,$ розглянемо величину
\begin{equation}\label{4}
{\cal E}_n(\mathfrak N)_X:=\sup\limits_{f\in \mathfrak N}\|f(\cdot)-{S}_{n-1}(f; \cdot)\|_X,
\end{equation}
в якій  $S_{n-1}(f; \cdot)$ --- частинна сума Фур'є порядку $n-1$ функції $f$.

При $X=C$ для точних верхніх меж вигляду (\ref{4}) на класах Вейля--Надя $W^r_{\beta,\infty}$  має місце асимптотична при $n\rightarrow\infty$ рівність 
\begin{equation}\label{4'}
{\cal E}_{n}(W^r_{\beta,\infty})_{C}=\frac4{\pi^2}\frac{\ln n}{n^rзнам}+O\left(\frac1{n^r}\right),\quad r>0,\quad\beta\in\mathbb R.
\end{equation}
При $r\in\mathbb{N}$ i $\beta=r$   цю оцінку довів А.М. Колмогоров \cite{Kolmogorov1985}. При дробових $r>0$ і деяких співвідношеннях між $r$ i $\beta$ --- В.Т.~Пінкевич \cite{Pinkevich1940} та С.М.~Нікольський \cite{Nikol'skii1941}. В загальному випадку оцінка (\ref{4'}) випливає із результатів А.В.~Єфімова \cite{Efimov1960} та С.О.~Теляковського \cite{Telyakovskii1961}. 

Зазначимо також, що аналогічна до (\ref{4'}) асимптотична рівність має місце і для класів $W^r_{\beta,1}$ в метриці простору $L_1$, а саме
\begin{equation}\label{5}
{\cal E}_{n}(W^r_{\beta,1})_{L_1}=\frac4{\pi^2}\frac{\ln n}{n^rзнам}+O\left(\frac1{n^r}\right),\quad r>0,\quad\beta\in\mathbb R,
\end{equation}
(див. \cite{Nikol'skii1946, Stechkin_Telyakovskii1967}).

У вказаних роботах параметри $r$ i $\beta$ класів Вейля--Надя вважалися фіксованими і питання про залежність залишкового члена в оцінці (\ref{4'}) (чи (\ref{5})) не розглядався.
Характер залежності від $r$ i $\beta$ залишкового члена в асимптотичній оцінці (\ref{4'}) вивчались  роботах І.Г. Соколова \cite{Sokolov1955}, С.Г. Селіванової \cite{Selivanova1955}, Г.І. Натансона \cite{Natanson1961}, С.О. Теляковського \cite{Telyakovskii1968, Telyakovskii1989, Telyakovskii1994} та С.Б. Стєчкіна \cite{Stechkin1980}.

У роботі \cite{Stechkin1980} досліджено асимптотичну поведінку величин ${\cal E}_{n}(W^r_{\beta,\infty})_{C}$ при $n\rightarrow\infty$ i $r\rightarrow\infty$. А саме було доведено, що при довільних $r\ge1$ i $\beta\in\mathbb{R}$ має місце рівність
\begin{equation}\label{6}
{\cal E}_{n}(W^r_{\beta,\infty})_{C}=\frac1{n^r}\left(\frac8{\pi^2}\mathbf{K}(e^{-r/n})+O(1)\frac1r\right),
\end{equation}
в якій
\begin{equation}\label{7}
\mathbf{K}(q)=\int\limits_{0}^{\pi/2}\frac{dt}{\sqrt{1-q^2\sin^2t}}
\end{equation}
--- повний еліптичний інтеграл першого роду, а $O(1)$ --- величина, рівномірно обмежена по $r, n$ i $\beta$.

Крім того, С.Б. Стєчкін показав \cite[теорема 4]{Stechkin1980}, що для великих $r$ залишковий член у рівності (\ref{6}) можна покращити. А саме, при довільних $r\ge n+1$ i $\beta\in\mathbb{R}$ виконується оцінка
\begin{equation}\label{8}
{\cal E}_{n}(W^r_{\beta,\infty})_{C}=
\frac1{n^r}\left(\frac4{\pi} +O(1)\left(1+\frac1n\right)^{-r} \right),
\end{equation}
в якій $O(1)$ --- величина, рівномірно обмежена по $r, n$ i $\beta$. Формула (\ref{8}) є асимптотичною рівністю, якщо $r/n\rightarrow\infty$.

З \cite{Stechkin1980} також випливає, що для величин ${\cal E}_{n}(W^r_{\beta,1})_{L_1}$ мають місце аналогічні до (\ref{6}) і (\ref{8}) оцінки, а саме при  $r\ge 1$ i $\beta\in\mathbb{R}$ рівномірно відносно всіх розглядуваних параметрів виконується формула
\begin{equation}\label{9}
{\cal E}_{n}(W^r_{\beta,1})_{L_1}=\frac1{n^r}\left(\frac8{\pi^2}\mathbf{K}(e^{-r/n})+O(1)\frac1r\right),
\end{equation}
а при $r\ge n+1$ i $\beta\in\mathbb{R}$ рівномірно по всіх параметрах --- оцінка
\begin{equation}\label{10}
{\cal E}_{n}(W^r_{\beta,1})_{L_1}=\frac1{n^r}\left(\frac4{\pi} +O(1)\left(1+\frac1n\right)^{-r} \right).
\end{equation}

Згодом С.О. Теляковський \cite{Telyakovskii1989,Telyakovskii1994} показав, що другий доданок у формулах \eqref{8} та \eqref{10}  можна замінити меншим, а саме, замість $\displaystyle O(1)(1+ 1/n)^{-r}$ записати $\displaystyle O(1)(1+ 2/n)^{-r}$. Він же уточнив оцінки \eqref{8} та \eqref{10} за рахунок виділення наступних членів асимптотики.

У роботі авторів \cite{Serdyuk_Sokolenko2019MFAT}, зокрема,  при довільних $1\le p\le\infty$ встановлено узагальнюючі аналоги оцінок  С.Б. Стєчкіна   \eqref{8} і \eqref{10}.
А саме, було розглянуто більш загальні, ніж $W^r_{\beta,p}$ класи функцій $W^r_{\bar\beta,p}$, які задаються згортками
\begin{equation}\label{11'}
	f(x)=\frac{a_0}{2}+\frac{1}{\pi}\int\limits_{-\pi}^{\pi}
	\varphi(x-t) B_{r,\bar\beta}(t)dt, \ \ \ a_0\in\mathbb R,\quad \varphi\in B^0_p,
\end{equation}
з ядрами 
\begin{equation}\label{11''}
	B_{r,\bar\beta}(t)=\sum\limits_{k=1}^\infty k^{-r}\cos\left(kt-\frac{\beta_k\pi}2\right),
\end{equation} 
що визначаються параметром $r>0$ і довільною числовою послідовністю $\bar\beta=\{\beta_k\}_{k=1}^\infty$ фазових зсувів $\beta_k\in\mathbb{R}$. 

Зрозуміло, що у випадку, коли $\beta_k=\beta$ при всіх $k\in\mathbb{N}$, то класи $W^r_{\bar\beta,p}$ перетворюються у класи Вейля--Надя $W^r_{\beta,p}$. 

У роботі \cite{Serdyuk_Sokolenko2019MFAT} для величин \eqref{4} при $\mathfrak{N}=W^r_{\bar\beta,p}$ i $X=C$ або $\mathfrak{N}=W^r_{\bar\beta,1}$ i $X=L_p$ було доведено, що при всіх 
 $r\ge n+1$, $\beta_k\in\mathbb{R}$ i  $1\le p\le\infty$ виконуються оцінки
\begin{equation}\label{11}
{\cal E}_{n}(W^r_{\bar\beta,p})_{C}=
\frac1{n^r}\left(\frac{\|\cos t\|_{p'}}{\pi} +O(1)\left(1+\frac1n\right)^{-r}\right),\quad\frac1p+\frac1{p'}=1,
\end{equation}
\begin{equation}\label{12}
{\cal E}_{n}(W^r_{\bar\beta,1})_{L_p}=
\frac1{n^r}\left(\frac{\|\cos t\|_{p}}{\pi} +O(1)\left(1+\frac1n\right)^{-r} \right),
\end{equation}
в яких $O(1)$ --- величини,	рівномірно обмежені відносно усіх розлядуваних параметрів. Оцінки (\ref{11}) і (\ref{12}) є асимптотичними рівностями, якщо $r/n\rightarrow\infty$.
Більш того, як доведено в \cite[теорема 5]{Serdyuk_Sokolenko2020UMB} за виконання умови $r/n\rightarrow\infty$ величини ${\cal E}_{n}(W^r_{\bar\beta,p})_{C}$ асимптотично збігаються з найкращими рівномірнми наближеннями класів $W^r_{\bar\beta,p}$ тригонометричними поліномами поряку $n-1.$

У роботі авторів \cite{Serdyuk_Sokolenko2019UMJ} встановлено інтерполяційні аналоги асимптотичних оцінок (\ref{11}) і (\ref{12}), в яких замість норм віхилень частинних сум Фур'є розглядаються норми відхилень інтерполяційних тригонометричних поліномів Лагранжа з рівномірним розподілом вузлів. 

В даній роботі досліджується асимптотична поведінка величин ${\cal E}_{n}(W^r_{\beta,p})_{C}$ i ${\cal E}_{n}(W^r_{\beta,1})_{L_p}, 1\le p\le\infty,$ у випадку, коли гладкісний та апроксимативний параметри $r$ i $n$ пов'язані співвідношенням
\begin{equation}\label{*}
	\sqrt{n}+1\le r\le n^2.
\end{equation} 

При $p=\infty$ асимптотична поведінка для зазначених величин відома і описується формулами \eqref{6} і \eqref{9}.

При $p=2$ відомі точні значення величин ${\cal E}_{n}(W^r_{\bar\beta,p})_{C}$ та ${\cal E}_{n}(W^r_{\bar\beta,1})_{L_p}$ при довільних $r>\frac12, n\in\mathbb{N}$ i $\bar\beta=\{\beta_k\}_{k=1}^\infty, \beta_k\in\mathbb{R}:$
\begin{equation}\label{romb}
	{\cal E}_{n}(W^r_{\bar\beta,2})_{C}={\cal E}_{n}(W^r_{\bar\beta,1})_{L_2}=\frac1{\sqrt\pi}\left(\sum\limits_{k=n}^\infty\frac1{k^{2r}}\right)^{1/2}=
	\frac1{\sqrt{\pi\Gamma(2r)}}\left(\int\limits_0^\infty\frac{t^{2r-1}e^{-nt}}{1-e^{-t}}dt\right)^{\frac12}
\end{equation}  
(див, наприклад, \cite{Serdyuk_Sokolenko2011,Serdyuk_Sokolenko2013,Serdyuk_Sokolenko2019MFAT}).

При всіх інших значеннях параметра $p$ (тобто при $1\le p<2$ або $2<p<\infty$ і $\frac rn\nrightarrow\infty$) асимптотичні рівності для величин ${\cal E}_{n}(W^r_{\beta,p})_{C}$ i ${\cal E}_{n}(W^r_{\beta,1})_{L_p}$ за виконання умови \eqref{*} залишались невідомими.

У даній роботі зазначені рівності будуть знайдені за рахунок використання методології, яка дозволяє звести задачу про сильну асимптотику величин вигляду \eqref{4} при $X=C$ або $L_p$ на класах Вейля--Надя до аналогічних величин на класах інтегралів Пуассона; для останніх асимптотичні рівності відомі завдяки роботам \cite{Serdyuk_2005_8,Serdyuk_2005_10}.

Позначимо через $C^{q}_{\beta,p},$ $q\in(0,1), \beta\in\mathbb{R}, 1\le p\le\infty, $ --- класи інтегралів Пуассона періодичних функцій з множин $B^0_p$ вигляду \eqref{3}, тобто класи $2\pi$-періодичних функцій $f$, що зображуються згортками 
\begin{equation}\label{1*}
	f(x)=\frac{a_0}{2}+\frac{1}{\pi}\int\limits_{-\pi}^{\pi}
	\varphi(x-t) P_{q,\beta}(t)dt, \quad a_0\in\mathbb R,\quad \varphi \in B_p^0,
\end{equation}
з ядрами Пуассона $P_{q,\beta}(t)$ вигляду
\begin{equation}\label{2*}
	P_{q,\beta}(t)=\sum\limits_{k=1}^\infty q^k \cos\left(kt-\frac{\beta\pi}2\right),\quad q\in(0,1),\quad  \beta\in\mathbb{R}.
\end{equation}

\textbf{2. Наближення функцій з класів \boldmath{$W^r_{\beta,p}$} сумами Фур'є в рівномірній метриці.}

\begin{theorem}\label{1t}
	Нехай $1 \le p \le \infty,$ $\beta\in\mathbb R,$  $r>1$ i $n\in\mathbb{N}$. Тоді за виконання умови 
	\begin{equation}\label{++}
		\sqrt{n}+1\le r\le n+1
	\end{equation} 
 при $p=1$ має місце формула
	\begin{equation}\label{1t1}
	{\cal E}_{n}(W^r_{\beta,1})_{C}=n^{-r}\Big(
	\frac1{\pi(1-e^{-r/n})}
	+O(1){n}{r^{-2}}\Big),	
	\end{equation}
	а при $1<p\le\infty$ 	--- формула
	\begin{equation}\label{1t2}
	{\cal E}_{n}(W^r_{\beta,p})_{C}=	
	n^{-r}\Big(
	\frac{\|\cos t\|_{p'}}{\pi}
	F^{\frac{1}{p'}}\Big(\frac{p'}{2},\frac{p'}{2};1;e^{-2r/n}\Big)+O(1){n}{r^{-2}}\Big),  \quad \frac1p+\frac1{p'}=1,
	\end{equation}
	де  $F(a,b;c;z)$ --- гіпергеометрична функція Гаусса
	\begin{equation}\label{1t3}
	F(a,b;c;z)=1+\sum\limits_{k=1}^{\infty}\frac{(a)_{k}(b)_{k}}{(c)_{k}}\frac{z^{k}}{k!},
	\ \ \ 
	(x)_{k}:=x(x+1)(x+2)...(x+k-1).
	\end{equation}
В \eqref{1t1} i \eqref{1t2} $O(1)$ --- величини, рівномірно обмежені відносно всіх розглядуваних параметрів.
\end{theorem}

\textit{\textbf{Доведення теореми \ref{1t}.}} Будемо використовувати запропонований С.Б.~Стєчкіним у \cite{Stechkin1980} метод доведеня, який полягає у тому, що залишок ряду Фур'є ядра Вейля--Надя $B_{r,\beta}$ вигляду \eqref{2} апроксимується в $L_{p'}$-метриці залишком ряду Фур'є ядра Пуассона $P_{q,\beta}$ вигляду \eqref{2*} при $q=e^{-r/n}$. Стєчкін \cite{Stechkin1980} реалізував зазначений підхід при $p'=1$. Ми ж розглядаємо загальний випадок $1 \le p'\le \infty.$ 

Покладемо 
\begin{equation}\label{1td1}
B_{r,\beta,n}(t):=\sum_{k=n}^{\infty}k^{-r}\cos\left(kt-\frac{\beta\pi}{2}\right).
\end{equation}
Із рівностей \eqref{1}, (\ref{2}) і (\ref{4}) отримуємо
\begin{equation}\label{1td2}
{\cal E}_{n}(W^r_{\beta,p})_{C}=\frac1{\pi}
\sup \limits_{\varphi\in B_p^0} \int\limits_{-\pi}^\pi
\varphi(x-t) \sum\limits_{k=n}^\infty k^{-r}\cos \left(nt-
\frac{\beta\pi}{2}\right)dt=\frac1{\pi} \sup\limits_{\varphi\in B_p^0}
\int\limits_{-\pi}^\pi \varphi(x-t)B_{r,\beta,n}(t)dt.
\end{equation}

Поклавши $q=e^{-r/n}$, подамо $B_{r,\beta,n}(t)$  у такому вигляді:
\begin{equation}\label{1td2'}
B_{r,\beta,n}(t)=\left(\frac en\right)^{r}P_{q,\beta,n}(t)+R_n(r;\beta)(t),
\end{equation}
де 
\begin{equation}\label{1td4}
P_{q,\beta,n}(t):=\sum_{k=n}^{\infty}q^k\cos\left(kt-\frac{\beta\pi}{2}\right),
\end{equation}
\begin{equation}\label{1td4**}
R_n(r;\beta)(t):=
B_{r,\beta,n}(t)-
\left(\frac en\right)^{r} P_{q,\beta,n}(t).
\end{equation}

Тоді з \eqref{1td2} одержуємо 
\begin{equation}\label{1td3}
{\cal E}_{n}(W^r_{\beta,p})_{C}=\Big(\frac en\Big)^{r}{\cal E}_{n}(C^q_{\beta,p})_{C}+O(1)R_n(r;\beta;p)_C,
\end{equation}
де
\begin{equation}\label{1td4'}
R_n(r;\beta;p)_C:=
\sup\limits_{\varphi\in B_p^0}
\int\limits_{-\pi}^\pi \varphi(x-t)R_n(r;\beta)(t)dt.
\end{equation}	

Як випливає з  \cite[теорема 1]{Serdyuk_2005_8} та \cite[формула (25)]{Serdyuk_2012_5} при довільних $1\leq p\leq\infty$ для величин ${\cal E}_{n}(C^q_{\beta,p})_{C}$, $0<q<1,$ $\beta\in\mathbb{R}$,
мають місце  рівності
\begin{equation}\label{1td5}
{\cal E}_{n}(C^q_{\beta,1})_{C}=q^n\Big(
\frac1{\pi(1-q)}
+O(1)\frac{q}{n(1-q)^2}\Big),	\quad p=1,
\end{equation}
та 	
\begin{equation}\label{1td5'}
{\cal E}_{n}(C^q_{\beta,p})_{C}=q^n\Big(
\frac{\|\cos t\|_{p'}}{\pi}
F^{\frac{1}{p'}}\Big(\frac{p'}{2},\frac{p'}{2};1;q^2\Big)+O(1)\frac{\xi(p)q}{n(1-q)^{s(p)}}\Big),  \quad 1<p\le\infty,  \frac1p+\frac1{p'}=1,
\end{equation}
де  
\begin{equation}\label{1td6}
s(p)=\left\{
\begin{array}{rl}
1, & p=\infty, \\
2, & p\in[1,\infty),
\end{array}
\right.
\quad 
\xi(p)={\left\{\begin{array}{cc}
	0, \ & p=2,  \ \ \ \ \ \ \ \ \ \ \ \ \ \ \ \ \  \\
	1, \ & p\in[1,2)\cup(2,\infty],
	\end{array} \right.}
\end{equation}
а величини $O(1)$ рівномірно обмежені відносно $n$, $p$, $q$ i $\beta$.

Для залишкових членів зі співвідношень (\ref{1td5}) і (\ref{1td5'}) при  $q=e^{-r/n}$ за умови $\sqrt{n}+1\le r\le n+1$, з урахуванням нерівностей
\begin{equation}\label{1td6'}
e^{-x}\le \frac1{1+x},\quad x>0,
\end{equation}
\begin{equation}\label{1td6**}
\frac1{1-e^{-x}}\le \frac{1+x}x,\quad x>0,
\end{equation}
отримуємо такі оцінки:  
\begin{equation}\label{1td6''}
\frac{q}{n(1-q)^{s(p)}}\le
\frac{e^{-r/n}}{n(1-e^{-r/n})^2}\le \frac1{n(1+r/n)}\left(\frac{1+r/n}{r/n}\right)^2=\frac{n+r}{r^2}\le\frac{2n+1}{r^2}=O(1){n}r^{-2}.	
\end{equation} 

Встановимо оцінку зверху залишку $R_n(r;\beta;p)_C$ зі співвідношення (\ref{1td3}). 
Покажемо, що за умови \eqref{++} для залишку  $R_n(r;\beta;p)_C$, означеного рівністю
\eqref{1td4'}, при всіх $1\le p\le \infty$ має місце рівномірна відносно усіх розглядуваних параметрів оцінка
\begin{equation}\label{1td6*}
	R_n(r;\beta;p)_C=O(1){n^{-r+1}}r^{-2}, \quad1\le p\le \infty, \quad\beta\in\mathbb{R}.
\end{equation}

Застосовуючи нерівність Гельдера до правої частини (\ref{1td4'}) та враховуючи, що $q=e^{-r/n}$, маємо
\begin{equation}\label{1td7}
R_n(r;\beta;p)_C\le \bigg\|B_{r,\beta,n}(t)-
\Big(\frac en\Big)^{r} P_{q,\beta,n}(t)\bigg\|_{p'}
=n^{-r}\left\|\sum_{k=n}^{\infty}\left(\left(\frac nk\right)^{r}-q^{k-n}\right)\cos\left(kt-\frac{\beta\pi}{2}\right)\right\|_{p'}=
$$
$$
=n^{-r}\left\|\sum_{k=1}^{\infty}\left(\left(1+\frac kn\right)^{-r}-e^{-rk/n}\right)\cos\left((k+n)t-\frac{\beta\pi}{2}\right)\right\|_{p'}\le
{(2\pi)^{1/p'}}n^{-r}
\sum_{k=1}^\infty \varphi\left(\frac kn\right),
\end{equation}
де
\begin{equation}\label{1td8}	
\varphi(x):=(1+x)^{-r}-e^{-rx}.
\end{equation}

Для оцінки зверху величини $\sum\limits_{k=1}^{\infty}\varphi\left(\frac kn\right)$ нам знадобиться твердження, доведення якого наведемо після доведення теореми 1.

\begin{lemma}\label{1l}
	Нехай  $r>1, n\in\mathbb{N}$ i    виконується умова \eqref{++}. 
Тоді мають місце нерівності 
	\begin{equation}\label{1l2}
	\sum_{k=1}^{\infty}\varphi\left(\frac kn\right)< (54e^{-1}+16-8\sqrt2)\,{n}r^{-2}<24,5518\,{n}r^{-2}. 
	\end{equation}
\end{lemma}

Оцінка \eqref{1td6*} у випадку виконання \eqref{++} є наслідком формул \eqref{1td7} і \eqref{1l2}.
Тоді, об'єднуючи співвідношення (\ref{1td3}), (\ref{1td5}), (\ref{1td5'}) i (\ref{1td6''}), отримуємо оцінки (\ref{1t1}) і (\ref{1t2}). 
Теорему \ref{1t} доведено.

\textit{\textbf{Доведення леми \ref{1l}.}} Покладемо $m= \left[\frac n{\sqrt{r}}\right],$ де $[x]$ --- ціла частина дійсного числа $x$. 
Зрозуміло, що 
\begin{equation}\label{d1l1}
m\le \frac n{\sqrt{r}}<m+1.
\end{equation}

Зазначимо, що (див. \cite[c. 145]{Stechkin1980})
\begin{equation}\label{d2l1}
0\le\varphi\left(\frac kn\right)\le e^{-\frac{rk}n}\left(e^{\frac{rk^2}{n(n+k)}}-1\right)=\left(e^{-\frac rn}\right)^k\left(e^{\frac{rk^2}{n(n+k)}}-1\right),\quad k\in\mathbb{N}.
\end{equation}

Крім того, згідно з теоремою Лагранжа, при $1\le k\le m$
\begin{equation}\label{d3l1}
\displaystyle e^{\frac{rk^2}{n(n+k)}}-1\le e^{\frac{rm^2}{n(n+1)}}\frac{rk^2}{n(n+k)}<e^{\frac{rm^2}{n^2}}\frac{rk^2}{n^2}.
\end{equation}

В силу (\ref{++}) $ {rm^2}/{n^2}\le 1$ і, отже, 
\begin{equation}\label{d4l1}
e^{\frac{rm^2}{n^2}}\le  e.
\end{equation}

Об'єднавши оцінки 	(\ref{d2l1})--(\ref{d4l1}) і врахувавши, що при довільних $0<q<1$
\begin{equation}\label{nabla}
	\sum_{k=1}^{\infty}q^kk^2=\frac{q(1+q)}{(1-q)^3},
\end{equation} 
одержуємо
\begin{equation}\label{d5l1}
\sum_{k=1}^{m}\varphi\left(\frac kn\right)\le\left(\sum_{k=1}^{m}\left(e^{-\frac rn}\right)^k k^2\right)\frac{er}{n^2}<\left(\sum_{k=1}^{\infty}\left(e^{-\frac rn}\right)^k k^2\right)\frac{er}{n^2}=
\frac{e^{-\frac rn}\left(1+e^{-\frac rn}\right)}{\left(1-e^{-\frac rn}\right)^3}\frac {er}{n^2}
<\frac{2e^{-\frac rn}}{\left(1-e^{-\frac rn}\right)^3}\frac {er}{n^2}.
\end{equation}

Оскільки, як неважко переконатись, функція $e^{-x}(1+x)^3$ приймає найбільше значення на $(0,+\infty)$ в точці $x_0=2$, тобто $\max\limits_{x>0}e^{-x}(1+x)^3=\frac{27}{e^2},$
то 
\begin{equation}\label{d6l1}
e^{-x}\le \frac{27}{e^2(1+x)^3}, \quad x>0.
\end{equation}

Скориставшись нерівностями \eqref{1td6**} i \eqref{d6l1}, маємо
\begin{equation}\label{d7l1}
\frac{e^{-\frac rn}}{\left(1-e^{-\frac rn}\right)^3}\le \frac{27}{e^2\left(1+\frac rn\right)^3}\frac{\left(1+\frac rn\right)^3}{\left(\frac rn\right)^3}=\frac{27}{e^2}\frac{n^3}{r^3}.
\end{equation}
Із (\ref{d5l1}) і (\ref{d7l1}) випливє нерівність 
\begin{equation}\label{d8l1}
\sum_{k=1}^{m}\varphi\left(\frac kn\right)<\frac{54}{e}\,{n}r^{-2}.
\end{equation}

Далі знайдемо оцінку зверху величини $\sum\limits_{k=m+1}^{\infty}\varphi\left(\frac kn\right).$ В силу  очевидної нерівності
\begin{equation}\label{d9l1}
\varphi\left(\frac kn\right)<\left(1+\frac kn\right)^{-r},
\end{equation}
та нерівності (\ref{d1l1}), маємо
\begin{equation}\label{d10l1'}
\varphi\left(\frac {m+1}n\right)<\left(1+\frac {m+1}n\right)^{-r}<\left(1+\frac{1}{\sqrt r}\right)^{-r}
\end{equation}
i
\begin{equation}\label{d10l1}
\sum_{k=m+2}^{\infty}\varphi\left(\frac kn\right)<\sum_{k=m+2}^{\infty}\left(1+\frac kn\right)^{-r}<\int\limits_{m+1}^{\infty}\left(1+\frac tn\right)^{-r}dt=n\int\limits_{\frac {m+1}n}^{\infty}\left(1+x\right)^{-r}dx=
$$
$$
=\frac n{r-1}\left(1+\frac{m+1}n\right)^{-r+1}<
\frac n{r-1}\left(1+\frac{1}{\sqrt r}\right)^{-r+1}.
\end{equation}
Функція $\frac{r^2}{r-1}\left(1+\frac1{\sqrt{r}}\right)^{-r+1}$ спадає на проміжку $[2,+\infty)$ і тому 
\begin{equation}\label{d11l1}
\max\limits_{r\ge2}\frac {r^2}{r-1}\left(1+\frac 1{\sqrt{r}}\right)^{-r+1}=8-4\sqrt2.
\end{equation}
Із (\ref{d10l1'}) i (\ref{d11l1}) отримуємо
\begin{equation}\label{d12l1}
\sum\limits_{k=m+1}^{\infty}\varphi\left(\frac kn\right)<\left(\frac{\sqrt r}{1+\sqrt{r}}+\frac n{r-1}\right) \left(1+\frac{1}{\sqrt r}\right)^{-r+1}= \frac{r+n-\sqrt r}{r-1}\left(1+\frac{1}{\sqrt r}\right)^{-r+1}< 
$$
$$
<\frac{2n+1-\sqrt r}{r-1}\left(1+\frac{1}{\sqrt r}\right)^{-r+1}<\frac {2n}{r^2}\cdot \frac{r^2}{r-1}\left(1+\frac 1{\sqrt{r}}\right)^{-r+1}\le (16-8\sqrt2)\,nr^{-2}.
\end{equation}

Об'єднавши оцінки 	(\ref{d8l1}) i (\ref{d12l1}), одержуємо нерівності (\ref{1l2}). Лему 1 доведено.

\begin{theorem}\label{2t}
	Нехай $1 \le p \le \infty,$ $\beta\in\mathbb R,$  $r>1$ i $n\in\mathbb{N}$. Тоді за виконання умови 
	\begin{equation}\label{+++}
		n+1\le r\le n^2
	\end{equation}  при $p=1$ має місце формула
	\begin{equation}\label{2t1}
	{\cal E}_{n}(W^r_{\beta,1})_{C}=n^{-r}\Big(
	\frac1{\pi(1-e^{-r/n})}
	+O(1)rn^{-2} e^{-r/n}\Big),	
	\end{equation}
	а при $1<p\le\infty$ --- формула
	\begin{equation}\label{2t2}
	{\cal E}_{n}(W^r_{\beta,p})_{C}=	
	n^{-r}\Big(
	\frac{\|\cos t\|_{p'}}{\pi}
	F^{\frac{1}{p'}}\Big(\frac{p'}{2},\frac{p'}{2};1;e^{-2r/n}\Big)+O(1)rn^{-2} e^{-r/n}\Big),  \quad \frac1p+\frac1{p'}=1,
	\end{equation}
	де  $F(a,b;c;z)$ --- гіпергеометрична функція Гаусса вигляду $(\ref{1t3})$. В \eqref{2t1} i \eqref{2t2} $O(1)$ --- величини, рівномірно обмежені відносно всіх розглядуваних параметрів.
\end{theorem}

\textit{\textbf{Доведення теореми \ref{2t}}} будемо проводити за схемою доведення теореми 1.  Базуючись на співвідношенні \eqref{1td3} і користуючись рівностями \eqref{1td5} і \eqref{1td5'}, неважко помітити, що асимптотичні формули \eqref{2t1} i \eqref{2t2} будуть встановлені, якщо за виконання \eqref{+++} доведемо істинність таких рівномірних по всіх параметрах оцінок:
\begin{equation}\label{48a}
	\frac{1}{(1-q)^2}=O(1), \quad \mbox{де}\quad q=e^{-r/n},
\end{equation}
та
\begin{equation}\label{48b}
	R_n(r;\beta;p)_C=O(1)rn^{-r-2} e^{-r/n}, \quad 1\le p\le\infty, \quad \beta\in\mathbb{R}.
\end{equation}
Щоб переконатись в справедливості \eqref{48a} досить скористатись нерівністю \eqref{1td6**} при $x=\frac rn$ та співвідношенням \eqref{+++} в силу яких
\begin{equation*} 
	\frac{1}{(1-q)^2}=\frac{1}{(1-e^{-r/n})^2}\le \left(\frac{1+r/n}{r/n}\right)^2=\left(\frac{r+n}{r}\right)^2\le\frac{2r-1}{r}<2.
\end{equation*}

Для оцінки залишку $R_n(r;\beta;p)_C$ використаємо ланцюжок співвідношень \eqref{1td7}, згідно з якими
\begin{equation}\label{49a}
	R_n(r;\beta;p)_C=O(1)\,n^{-r}\sum_{k=1}^\infty \varphi\left(\frac kn\right).
\end{equation}

Як випливає з доведення теореми 3 роботи С.Б.~Стєчкіна \cite[c.~147]{Stechkin1980} при $n\in\mathbb{N}$ i $ n+1\le r\le n^2$  рівномірно відносно усіх розглядуваних параметрів справедлива оцінка
\begin{equation}\label{2td2}
\sum_{k=1}^\infty \varphi\left(\frac kn\right)=O(1)rn^{-2} e^{-r/n}.
\end{equation}		

Оцінка \eqref{48b} випливає безпосередньо із формул \eqref{49a} та \eqref{2td2}.
Теорему \ref{2t} доведено.

\textbf{3. Наближення функцій з класів \boldmath{$W^r_{\beta,1}$} сумами Фур'є в інтегральних метриках.}

\begin{theorem}\label{3t}
	Нехай $1 \le p \le \infty,$ $\beta\in\mathbb R,$  $r>1$ i $n\in\mathbb{N}$. Тоді за виконання умови $\sqrt{n}+1\le r\le n+1$  при $1\le p<\infty$ має місце формула
	\begin{equation}\label{3t1}
	{\cal E}_{n}(W^r_{\beta,1})_{L_p}=	
	n^{-r}\Big(
	\frac{\|\cos t\|_{p}}{\pi}
	F^{\frac{1}{p}}\Big(\frac{p}{2},\frac{p}{2};1;e^{-2r/n}\Big)+O(1){n}{r^{-2}}\Big),  \quad \frac1p+\frac1{p'}=1,	
	\end{equation}
	а при  	$p=\infty$ --- формула
	\begin{equation}\label{3t2}
	{\cal E}_{n}(W^r_{\beta,1})_{L_\infty}=n^{-r}\Big(
	\frac1{\pi(1-e^{-r/n})}
	+O(1){n}{r^{-2}}\Big),
	\end{equation}
	де $F(a,b;c;z)$ --- гіпергеометрична функція Гаусса вигляду $(\ref{1t3})$,
	а $O(1)$ --- величини, рівномірно обмежені відносно всіх розглядуваних параметрів.
\end{theorem}

\textit{\textbf{Доведення теореми \ref{3t}.}} Аналогічно до того, як це реалізовано при доведенні теореми 1, подамо величину ${\cal E}_{n}(W^r_{\beta,1})_{L_p}, r>1, 1\le p\le\infty, \beta\in\mathbb{R},$ у вигляді
\begin{equation}\label{3td3}
{\cal E}_{n}(W^r_{\beta,1})_{L_p}=\Big(\frac en\Big)^{r}{\cal E}_{n}(C^q_{\beta,1})_{L_p}+O(1)R_n(r;\beta;1)_{L_p},
\end{equation}
де
\begin{equation}\label{3td4}
{\cal E}_{n}(C^q_{\beta,1})_{L_p}=\frac1{\pi} \sup\limits_{\varphi\in B_1^0}
\left\|\int\limits_{-\pi}^\pi \varphi(x-t)P_{q,\beta,n}(t)dt\right\|_p,\quad q=e^{-r/n},
\end{equation}	
\begin{equation}\label{3td4'}
R_n(r;\beta;1)_{L_p}:=
\sup\limits_{\varphi\in B_1^0}
\left\|\int\limits_{-\pi}^\pi \varphi(x-t)R_n(r;\beta)(t)dt\right\|_p,
\end{equation}
а $R_n(r;\beta)(t)$ означено рівністю \eqref{1td4**}.

Як випливає з теореми 1 роботи \cite{Serdyuk_2005_10} та формули (25) роботи \cite{Serdyuk_2012_5} при довільних $1\leq p\leq\infty$ для величин ${\cal E}_{n}(C^q_{\beta,1})_{L_p}$, $0<q<1,$ $\beta\in\mathbb{R}$,
мають місце асимптотичні рівності    
\begin{equation}\label{3td1}
{\cal E}_{n}(C^q_{\beta,1})_{L_p}=
\left\{\begin{array}{ll}
\displaystyle	q^n\bigg(\frac{\|\cos t\|_{p}}{\pi}F^{\frac{1}{p}}\Big(\frac{p}{2},\frac{p}{2};1;q^2\Big) +O(1)\frac{q}{n(1-q)^{s(p')}}\bigg), \ & 1\leq p<\infty,  \displaystyle\frac1p+\frac1{p'}=1,
	\\ \ & \ 
	\\
\displaystyle	q^n\bigg(\frac{1}{\pi(1-q)}+O(1)\frac{q}{n(1-q)^{2}}\bigg),  & p=\infty,
	\end{array} \right.
\end{equation}
в яких $s(\cdot)$ означається співвідношенням (\ref{1td6}), a $O(1)$ --- величини, що рівномірно обмежені відносно усіх розглядуваних параметрів.

В роботі \cite[с. 250]{Serdyuk_Sokolenko2013} було доведено, що для величин ${\cal E}_{n}(C^q_{\beta,1})_{L_p}$ при $p=2$  виконується рівність
\begin{equation}\label{3td1'}
{\cal E}_{n}(C^q_{\beta,1})_{L_2}=\frac{q^{n}}{\sqrt{\pi(1-q^2)}}, \quad 0<q<1, \ \beta\in\mathbb{R}, \ n\in\mathbb{N}.
\end{equation}

\noindent Рівність (\ref{3td1'}) уточнює асимптотичну рівність (\ref{3td1}) при $p=2$  в тому сенсі, що зазначена рівність (\ref{3td1}) при $s=2$ залишається вірною, якщо в ній обнулити залишковий член. 

Отже, взявши до уваги формули (\ref{3td1}), (\ref{3td1'})  та очевидну рівність $F(1,1;1;q^{2})=\displaystyle\frac{1}{1-q^{2}}, \ \ q\in(0,1)$, для всіх $0<q<1$, $\beta\in\mathbb{R}$ і $1\leq p\leq\infty$ можемо записати 
\begin{equation}\label{3td1''}
{\cal E}_{n}(C^q_{\beta,1})_{L_p}=
\left\{\begin{array}{lc}
\displaystyle q^n\bigg(\frac{\|\cos t\|_{p}}{\pi}F^{\frac{1}{p}}\Big(\frac{p}{2},\frac{p}{2};1;q^2\Big) +O(1)\frac{\xi(p)q}{n(1-q)^{s(p')}}\bigg), \ & 1\leq p<\infty,  \displaystyle\frac1p+\frac1{p'}=1,
\\ \ & \ 
\\
\displaystyle q^n\bigg(\frac{1}{\pi(1-q)}+O(1)\frac{q}{n(1-q)^{2}}\bigg),  & p=\infty,
\end{array} \right.
\end{equation}
в яких $\xi(\cdot)$ i $s(\cdot)$ означено формулами (\ref{1td6}), a $O(1)$ --- величини, що рівномірно обмежені відносно усіх розглядуваних параметрів.

Застосовуючи до правої частини (\ref{3td4'}) твердження 1.5.5 із роботи \cite[с.~43]{Kornejchuk1987}, згідно з яким для $L_p$-норми згортки $(\varphi*K)(\cdot)=\frac1\pi\int\limits_{-\pi }^\pi\varphi(\cdot-t)K(t)dt,$ де $\varphi\in L_1, K\in L_p,  1\le p\le \infty,$
\begin{equation}\label{3td5}
\|\varphi*K\|_p \leq \frac{1}{\pi} \|\varphi\|_1\|K\|_p,  
\end{equation} 
та враховуючи оцінку (\ref{1l2}), для  залишку $R_n(r;\beta;1)_{L_p}$ при $1\le p\le\infty, \beta\in\mathbb{R}$ і виконанні умови \eqref{++}, одержуємо 
\begin{equation}\label{3td7}
R_n(r;\beta;1)_{L_p}\le \left\|R_n(r;\beta)(t)\right\|_{p}
\le
{(2\pi)^{1/p}}n^{-r}
\sum_{k=1}^\infty \varphi\left(\frac kn\right)=O(1)n^{1-r}r^{-2}.
\end{equation}

Із  (\ref{1td6''}), (\ref{3td3}), (\ref{3td1''}) та (\ref{3td7})   отримуємо рівності (\ref{3t1}) і (\ref{3t2}).  Теорему \ref{3t} доведено.

\begin{theorem}\label{4t}
	Нехай $1 \le p \le \infty,$ $\beta\in\mathbb R,$  $r>1$ i $n\in\mathbb{N}$. Тоді за виконання умови $n+1\le r\le n^2$  при $1\le p<\infty$ має місце рівність
	\begin{equation}\label{4t1}
	{\cal E}_{n}(W^r_{\beta,1})_{L_p}=	
	n^{-r}\Big(
	\frac{\|\cos t\|_{p}}{\pi}
	F^{\frac{1}{p}}\Big(\frac{p}{2},\frac{p}{2};1;e^{-2r/n}\Big)+O(1)rn^{-2} e^{-r/n}\Big),	\quad \frac1p+\frac1{p'}=1,	
	\end{equation}
	а при $p=\infty$ --- рівність
	\begin{equation}\label{4t2}
	{\cal E}_{n}(W^r_{\beta,1})_{L_\infty}=n^{-r}\Big(
	\frac1{\pi(1-e^{-r/n})}
	+O(1)rn^{-2} e^{-r/n}\Big),   
	\end{equation}
	де  $F(a,b;c;z)$ --- гіпергеометрична функція Гаусса вигляду $(\ref{1t3})$,
	а $O(1)$ --- величини, рівномірно обмежені відносно всіх розглядуваних параметрів.
\end{theorem}

\noindent\textit{\textbf{Доведення теореми \ref{4t}.}}  За умови $n+1\le r\le n^2$ в силу \eqref{2td2}, \eqref{3td4'} і \eqref{3td5}
\begin{equation}\label{4td2}
R_n(r;\beta;1)_{L_p}=\|R_n(r;\beta)(t)\|_p\le 
{(2\pi)^{1/p}}n^{-r}
\sum_{k=1}^\infty \varphi\left(\frac kn\right)=O(1)rn^{-r-2} e^{-r/n}. 
\end{equation}

Об'єднуючи співвідношення  (\ref{3td3}), (\ref{3td1''}) та (\ref{4td2}), для величин ${\cal E}_{n}(W^r_{\beta,1})_{L_p}$ за виконання умови \eqref{+++} отримуємо рівності
\begin{equation*}\label{4td3}
{\cal E}_{n}(W^r_{\beta,1})_{L_p}=\left\{\begin{array}{lc}
\displaystyle n^{-r}\bigg(\frac{\|\cos t\|_{p}}{\pi}F^{\frac{1}{p}}\Big(\frac{p}{2},\frac{p}{2};1;e^{-2r/n}\Big) +O(1)rn^{-2} e^{-r/n}\bigg), \ & 1\leq p<\infty,  
\\ \ & \ 
\\
\displaystyle n^{-r}\bigg(\frac{1}{\pi(1-e^{-r/n})}+O(1)rn^{-2} e^{-r/n}\bigg),  & p=\infty,
\end{array} \right.
\end{equation*}
Теорему \ref{4t} доведено.

Зрозуміло, що у формулах \eqref{3t2} i \eqref{4t2} величини ${\cal E}_{n}(W^r_{\beta,1})_{L_\infty}$ можна замінити на ${\cal E}_{n}(W^r_{\beta,1})_C$
 і, по суті, зазначені оцінки збігаються з оцінками \eqref{1t1} i \eqref{2t1}, відповідно. Доцільність їх наведення у теоремах \ref{3t} та \ref{4t} мотивується завершеністю формулювання останніх відносно параметра $p\ (1\le p\le\infty)$.

В ході доведення теорем \ref{1t} -- \ref{4t} суттєвим чином були використані асимптотичні рівності для точних верхніх меж відхилень сум Фур'є на класах інтегралів Пуассона $C^q_{\beta,p}$ (формули \eqref{1td5}, \eqref{1td5'}, \eqref{3td1}). Завдяки роботам \cite{Stepanets1995,Stepanets2005,Telyakovskii1989,Telyakovskii1994,Serdyuk_Stepanyuk_2015UMJ,Serdyuk_Stepanyuk2019Anal} аналогічні асимптотичні рівності встановлено і для класів узагальнених інтегралів Пуассона $C^{\alpha,r}_{\beta,p}$ (див., наприклад, \cite{Serdyuk_Stepanyuk2019Anal}). 

\textbf{4. Зауваження та наслідки.}

\begin{remark}\label{1r}
Формули \eqref{1t1}, \eqref{1t2}, \eqref{2t1}, \eqref{2t2}, \eqref{3t1}, \eqref{3t2}, \eqref{4t1} i \eqref{4t2}, які фігурують в теоремах \ref{1t}--\ref{4t} є асимтотичними рівностями при $r\rightarrow\infty, n\rightarrow\infty$ у випадку, коли частка $\frac rn$ обмежена зверху і знизу деякими додатними числами $K_1$ i $K_2$:
\begin{equation}\label{4.1}
	0<K_1\le\frac rn\le K_2<+\infty.
\end{equation}	
\end{remark} 

Дійсно, нехай спочатку $\sqrt{n}+1\le r\le n+1$. В цьому випадку 
\begin{equation}\label{4.2}
		0<K_1\le\frac rn\le 2
\end{equation}
і, отже, можна записати
\begin{equation}\label{4.3}
	\frac n{r^2}=O\left(\frac1r\right).
\end{equation}

Далі, з урахуванням \eqref{4.2}, маємо
\begin{equation}\label{4.4}
	\frac1{1-e^{-r/n}}\ge\frac1{1-e^{-2}},
\end{equation}
\begin{equation}\label{4.5}
	F\Big(s,s;1;e^{-2r/n}\Big)\ge F\Big(s,s;1;e^{-4}\Big),\quad s>0.
\end{equation}
Із співвідношень \eqref{4.2}, \eqref{4.4} i \eqref{4.5} випливає, що за умови \eqref{4.1} формули \eqref{1t1}, \eqref{1t2}, \eqref{3t1} i \eqref{3t2}
є асимптотичними рівностями при $r\rightarrow\infty$ i $ n\rightarrow\infty.$ При цьому у \eqref{1t1} і \eqref{1t2} та \eqref{3t1} i \eqref{3t2} залишковий член $O(1)nr^{-2}$ можна замінити на $O(\frac1r)$.

Нехай, далі, $n+1\le r\le n^2$. В цьому випадку 
\begin{equation}\label{4.6}
	1<\frac rn\le K_2<+\infty,
\end{equation}
і, отже,
\begin{equation}\label{4.7}
	\frac r{n^2}e^{-r/n}=O\left(\frac1r\right).
\end{equation}
Очевидно також, що 
\begin{equation}\label{4.8}
\frac1{1-e^{-r/n}}>1,
\end{equation}
\begin{equation}\label{4.9}
F\Big(s,s;1;e^{-2r/n}\Big)>1,\quad s>0.
\end{equation}

Із співвідношень \eqref{4.7}, \eqref{4.8} i \eqref{4.9} випливає, що  формули \eqref{2t1}, \eqref{2t2}, \eqref{4t1} i \eqref{4t2} за умови \eqref{4.1} також 
є асимптотичними рівностями при $r\rightarrow\infty$ i $ n\rightarrow\infty.$ 

\begin{remark}\label{2r}
	Формули \eqref{1t1}, \eqref{1t2},  \eqref{3t1} i \eqref{3t2} з теорем \ref{1t} і \ref{3t} є асимтотичними рівностями при $r\rightarrow\infty, n\rightarrow\infty$ у випадку, коли 
	\begin{equation}\label{4.10}
		\lim\limits_{^{r\rightarrow\infty}_{n\rightarrow\infty}}\frac rn=0.
	\end{equation}	
\end{remark} 

Дійсно, в силу умови $\sqrt{n}+1\le r\le n+1$, що фігурує в зазначених теоремах, має місце нерівність 
\begin{equation}\label{4.11}
	\frac n{r^2}\le1,
\end{equation}
а, отже, залишкові члени у формулах \eqref{1t1}, \eqref{1t2},  \eqref{3t1} i \eqref{3t2} рівномірно обмежені по усіх параметрах. Залишається показати, що головні члени у зазначених формулах за виконання умови \eqref{4.10} прямують до нескінченності при $r\rightarrow\infty$ i $ n\rightarrow\infty.$ 
Оскільки 
\begin{equation}\label{4.12}
	\frac{1}{1-e^{-x}}=\frac1x+O(1)\quad \mbox{при}\quad x\rightarrow 0,
\end{equation}
то за умови \eqref{4.10} справджується асимптотична рівність 
\begin{equation}\label{4.13}
	\frac1{1-e^{-r/n}}=\frac nr+O(1).
\end{equation}

Із \eqref{4.10}, \eqref{4.11} i \eqref{4.13} випливає, що формули \eqref{1t1} i \eqref{3t1} є асимптотичними рівностями. Більш того, насправді ми довели наступне твердження.
\begin{corollary}\label{1c}
	Нехай $\sqrt n+1\le r\le n+1, n\in\mathbb{N}, \beta\in\mathbb{R}$ і має місце \eqref{4.10}. Тоді при $r\rightarrow\infty$ i $ n\rightarrow\infty$ виконується асимптотична рівність
	\begin{equation}\label{1c1}
		{\cal E}_{n}(W^r_{\beta,1})_C={\cal E}_{n}(W^r_{\beta,1})_{L_\infty}=\frac1{rn^{r-1}}
		\left(	\frac1\pi
		+O\left(\frac rn+\frac 1r\right)\right). 
	\end{equation}
\end{corollary}

Покажемо далі, що за умови \eqref{4.10}
\begin{equation}\label{4.14}
	\lim\limits_{^{r\rightarrow\infty}_{n\rightarrow\infty}}F^{1/s}\Big(\frac s2,\frac s2;1;e^{-2r/n}\Big)=+\infty,\quad s>0.
\end{equation}

Як показано в \cite[формула (25)]{Serdyuk_2012_5}
\begin{equation}\label{4.15}
	F^{1/s}\Big(\frac s2,\frac s2;1;q^2\Big)=\left(\frac1{2\pi}\int\limits_{0}^{2\pi}\left|\frac1{\sqrt{1-2q\cos x+q^2}}\right|^sdx\right)^{1/s},\quad s\ge1,\quad q\in(0,1),
\end{equation}
тому, враховуючи, що величина
$$
\left(\frac1{2\pi}\int\limits_{0}^{2\pi}\left|\frac1{\sqrt{1-2q\cos x+q^2}}\right|^sdx\right)^{1/s}
$$
зростає по параметру $s$ на $[1,+\infty)$, отримуємо, що
\begin{equation}\label{4.16}
	F^{1/s}\Big(\frac s2,\frac s2;1;q^2\Big)\ge F\Big(\frac 12,\frac 12;1;q^2\Big)=
	\frac1{2\pi}\int\limits_{0}^{2\pi}\frac{dx}{\sqrt{1-2q\cos x+q^2}}=\frac2\pi\mathbf{K}(q),
\end{equation}
де $\mathbf{K}(q)$ --- повний еліптичний інтеграл першого роду.

Із \eqref{4.16} і асимптотичного розкладу величини $\mathbf{K}(q)$
$$
\mathbf{K}(q)=\frac12\ln \frac1{1-q}+C+o(1),\quad q\rightarrow1-0
$$
(див. \cite[Гл.~22]{30}), а також формули \eqref{4.13} маємо
\begin{equation}\label{4.17}
	F^{1/s}\Big(\frac s2,\frac s2;1;e^{-2r/n}\Big)\ge \frac1\pi\ln \frac1{1-e^{-r/n}}+O(1)=\frac1\pi\ln\frac nr+O(1)\rightarrow\infty.
\end{equation}

Із \eqref{4.11} і \eqref{4.17} випливає, що за умови \eqref{4.10} оцінки \eqref{1t1} i \eqref{3t1} є асимптотичними рівностями при $r\rightarrow\infty$ i $ n\rightarrow\infty$.

\begin{remark}\label{3r}
	Формули \eqref{2t1}, \eqref{2t2},  \eqref{4t1} i \eqref{4t2}, які фігурують в теоремах \ref{2t} і \ref{4t} є асимтотичними рівностями при $r\rightarrow\infty, n\rightarrow\infty$ у випадку, коли 
	\begin{equation}\label{4.18}
		\lim\limits_{^{r\rightarrow\infty}_{n\rightarrow\infty}}\frac rn=+\infty.
	\end{equation}	
\end{remark} 

Дійсно, як зазначалось раніше, в силу умови $n+1\le r\le n^2$ має місце оцінка \eqref{4.7}  і, крім того, в силу \eqref{4.18}
\begin{equation}\label{4.19}
	\frac1{1-e^{-r/n}}=1+e^{-r/n}+O(1)e^{-2r/n}.
\end{equation}

Із \eqref{4.7} і \eqref{4.19} випливає, що формули   \eqref{2t1}  i \eqref{4t2} є асимптотичними рівностями. Більш того, ми тим самим довели наступне твердження.

\begin{corollary}\label{2c}
	Нехай $ n+1\le r\le n^2, n\in\mathbb{N}, \beta\in\mathbb{R}$ і має місце \eqref{4.18}. Тоді при $r\rightarrow\infty$ i $ n\rightarrow\infty$ виконується асимптотична рівність
	\begin{equation}\label{2c1}
		{\cal E}_{n}(W^r_{\beta,1})_C={\cal E}_{n}(W^r_{\beta,1})_{L_\infty}=\frac1{n^r}		\left(	\frac1\pi+\frac1\pi e^{-r/n}
		+O(1)\left(\frac r{n^2}+e^{-r/n}\right) e^{-r/n}\right). 
	\end{equation}
\end{corollary}

Із співвідношень \eqref{4.15} і \eqref{4.16} для величини $F^{1/s}\Big(\frac s2,\frac s2;1;q^2\Big)$ випливає двостороння оцінка
\begin{equation}\label{4.20}
\frac2\pi\mathbf{K}(q)\le F^{1/s}\Big(\frac s2,\frac s2;1;q^2\Big)\le\frac1{1-q},\quad s\ge1,\quad q\in(0,1).
\end{equation}

Оскільки в силу \eqref{4.16}
\begin{equation}\label{4.21}
	\mathbf{K}(q)=\frac\pi2+O(q^2),
\end{equation}
то із \eqref{4.20} при $q=e^{-r/n}$ за умови \eqref{4.18} випливає асимптотична при 
$r\rightarrow\infty$ i $ n\rightarrow\infty$ рівність
\begin{equation}\label{4.22}
	F^{1/s}\Big(\frac s2,\frac s2;1;e^{-2r/n}\Big)=1+O(1)e^{-r/n},\quad 1\le s<\infty.
\end{equation}

Із \eqref{4.7}  і \eqref{4.22} випливає, що за умови \eqref{4.18}  
формули   \eqref{2t2}  i \eqref{4t1} є асимптотичними при $r\rightarrow\infty$ i $ n\rightarrow\infty$ рівностями. При цьому нами встановлено наступне твердження.

\begin{corollary}\label{3c}
	Нехай $1\le p\le\infty,$ $n\in\mathbb{N}, \beta\in\mathbb{R}$ i $ n+1\le r\le n^2$. Тоді за виконання умови умови \eqref{4.18} мають місце асимптотичні  при $r\rightarrow\infty$ i $ n\rightarrow\infty$ рівності
\begin{equation}\label{3c1}
	{\cal E}_{n}(W^r_{\beta,p})_{C}=
	\frac1{n^r}\left(\frac{\|\cos t\|_{p'}}{\pi} +O(1)e^{-r/n}\right),\quad\frac1p+\frac1{p'}=1,
\end{equation}
\begin{equation}\label{3c2}
	{\cal E}_{n}(W^r_{\beta,1})_{L_p}=
	\frac1{n^r}\left(\frac{\|\cos t\|_{p}}{\pi} +O(1)e^{-r/n} \right),
\end{equation}
\end{corollary}

Утім асимптотичні рівності \eqref{3c1}  і \eqref{3c2} випливають також із \eqref{11}  і \eqref{12}, відповідно. Отже формули \eqref{11}  і \eqref{12} та  \eqref{2t1}, \eqref{2t2},  \eqref{4t1}  i \eqref{4t2}  повністю узгоджуються між собою.

\begin{remark}\label{4r}
	Оцінки \eqref{1t1} i \eqref{2t1} є граничними випадками оцінок \eqref{1t2} i \eqref{2t2} при $p'\rightarrow\infty$. Аналогічно, оцінки \eqref{3t2} i \eqref{4t2} є граничними випадками оцінок \eqref{3t1} i \eqref{4t1} при $p\rightarrow\infty$.
\end{remark} 

Щоб у цьому переконатись треба перейти до границі при $s\rightarrow\infty$ у формулі 
\eqref{4.15}:
\begin{equation}\label{4.23}
\lim_{s\rightarrow\infty}	F^{1/s}\Big(\frac s2,\frac s2;1;q^2\Big)=\lim_{s\rightarrow\infty}\frac1{(2\pi)^{1/s}}\left(\int\limits_{0}^{2\pi}\left|\frac1{\sqrt{1-2q\cos x+q^2}}\right|^sdx\right)^{1/s}=
$$
$$
=
\left\|\frac1{\sqrt{1-2q\cos (\cdot)+q^2}}\right\|_{L_\infty}=\frac1{1-q}.
\end{equation}

Далі залишається застосувати \eqref{4.23} при $q=e^{-r/n}$ i $s=p'$ або $s=p$.

\small{
	
}

\end{document}